\documentclass{article}
\usepackage{amssymb,xy}

\xyoption{arrow}
\xyoption{curve}
\xyoption{frame}
\xyoption{matrix}

\newtheorem{theorem}{Theorem}
\newtheorem{proposition}[theorem]{Proposition}
\newtheorem{lemma}[theorem]{Lemma}
\newtheorem{corollary}[theorem]{Corollary}
\newtheorem{remark}[theorem]{Remark}

\def\proof{\smallskip\noindent {\it Proof --- \ }}
\def\proofof#1{\smallskip\noindent {\it Proof of #1 --- \ }}
\def\endproof{\hfill$\square$\medskip}

\def\C{\mathbb{C}}
\def\H{\mathrm{H}} 
\def\l{\ell}
\def\QH{\mathrm{QH}} 
\def\bQH{\overline{\mathrm{QH}}\kern0.01cm}
\def\qp{\,*\,}

\def\s{\sigma}
\def\wnot{{w_{\mathrm{o}}}}
\def\Z{\mathbb{Z}}
\def\<{\left<}
\def\>{\right>}

\title{Symmetries of Gromov-Witten Invariants}

\author{Alexander Postnikov \\[.05in]
{\small Department of Mathematics,
University of California, Berkeley, CA 94720}\\
{\small \texttt{apost@math.berkeley.edu}}}

\date{September 17, 2000}

\begin{document}
\maketitle

\begin{abstract} 
The group $(\Z/n\Z)^2$ is shown to act on the Gromov-Witten 
invariants of the complex flag manifold.  We also deduce 
several corollaries of this result.
\end{abstract}

\section{Introduction}
\label{se:intro}

The aim of this paper is to present certain symmetry properties of
the Gromov-Witten invariants for type~$A$ complex flag manifolds.

Recall that the cohomology ring of the complex flag manifold $Fl_n$
has an additive basis of Schubert classes $\s_w$, which are indexed 
by permutations $w$ in the symmetric group $S_n$.
For permutations $u,v,w\in S_n$, the Schubert number~$c_{u,v,w}$
is the structure constant of the cohomology ring 
in the basis of Schubert classes:
$$
\s_u\cdot \s_v = \sum_{w\in S_n} c_{u,v,w} \,\s_{\wnot w}\,,
$$
where $\wnot$ is the longest permutation in $S_n$.  Equivalently,
$$
c_{u,v,w} = \int \s_u \cdot \s_u \cdot \s_w
$$
is the intersection number of Schubert varieties.
Thus these numbers are nonnegative integers symmetric in~$u$, $v$, and~$w$.
They generalize the famous Littlewood-Richardson coefficients.
If $\l(u)+\l(v)+\l(w)\ne {n(n-1)\over 2}$ then the Schubert number
$c_{u,v,w}$ is zero for an obvious degree reason.

A long standing open problem is to find an algebraic or combinatorial
construction for the coefficients $c_{u,v,w}$ that would imply their
nonnegativity.  A possible approach to this problem could be in its
generalization to the quantum cohomology ring of the flag manifold $Fl_n$.
The structure constants of this ring are certain polynomials whose 
coefficients are the Gromov-Witten invariants
$\<\s_u,\s_v,\s_w\>_{(d_1,\dots,d_{n-1})}$
The Schubert number $c_{u,v,w}$ is
a special case of the Gromov-Witten invariants: 
$c_{u,v,w}=\<\s_u,\s_v,\s_w\>_{(0,\dots,0)}$.   
These invariants are defined as numbers of certain rational curves 
in $Fl_n$.  The geometric definition of the
Gromov-Witten invariants implies their nonnegativity.

In this paper we establish cyclic symmetries of the Gromov-Witten invariants
that could not be detected in their full generality on the ``classical'' level 
of the Schubert numbers $c_{u,v,w}$.
Several related results for the~$c_{u,v,w}$ when $u$ is a Grassmannian
permutation were, however,
found by Bergeron and Sottile, see~\cite[Theorems~1.3.4, 1.3.4]{BS}.  In case
of the Gromov-Witten invariants we do not need to restrict the rule to
Grassmannian permutations.  Similar symmetries of the Gromov-Witten invariants 
for Grassmannian varieties were found in~\cite{AW}.

\section{Gromov-Witten invariants}
\label{se:GW}  

Let~$Fl_n$ denote the manifold of complete flags of subspaces in 
the complex $n$-dimensional linear space~$\C^n$.  One can also define
the {\it flag manifold\/} as $Fl_n=GL_n(\C)/B$, where $B$ is the Borel
subgroup of upper triangular matrices in the general linear group.
The flag manifold is a compact smooth complex manifold.
For a permutation $w\in S_n$, the {\it Schubert variety\/} $X_w$ 
is the closure of the {\it Schubert cell\/} $B_{-}wB/B$ in $Fl_n$, where
$B_{-}$ is the subgroup of lower triangular matrices and $w$ is viewed
as a permutation matrix in $GL_n$.
The {\it Schubert classes\/}~$\s_w\in \H^*(Fl_n,\Z)$, indexed by 
permutations~$w\in S_n$, are defined as the Poincar\'e duals of the homology
classes $[X_w]$ of Schubert manifolds.  They form an additive 
$\Z$-basis of the cohomology ring~$\H^*(Fl_n,\Z)$.  
Moreover, $\s_w\in \H^{2l}(Fl_n,\Z)$, where $l=\l(w)$
is the {\it length} of permutation~$w$, i.e., its number of inversions.

Recently, attention has been drawn to the (small) {\it quantum cohomology
ring\/} $\QH^*(Fl_n,\Z)$ of the flag manifold.  The definition of quantum 
cohomology can be found, for example, in~\cite{FP}.   
Here we briefly outline several notions and results.

As a vector space, the quantum cohomology of~$Fl_n$
is the usual cohomology tensored with the polynomial ring in $n-1$ 
variables:
\begin{equation}
\QH^*(Fl_n,\Z)\cong \H^*(Fl_n,\Z)\otimes \Z[q_1,\dots,q_{n-1}].
\label{eq:QH}
\end{equation}
The Schubert classes $\s_w$, thus, form a $\Z[q_1,\dots,q_{n-1}]$-basis 
of the quantum cohomology ring.

The multiplication in~$\QH^*(Fl_n,\Z)$ (quantum product) is a commutative
$\Z[q_1,\dots,q_{n-1}]$-linear operation.   It is therefore sufficient to specify the quantum
product of any two Schubert classes.  To avoid confusion with the
multiplication in the usual cohomology ring, we will use ``$*$'' to denote 
the quantum product.
The quantum product $\s_u \qp \s_v$ of two Schubert classes 
can be expressed  in the basis of the Schubert classes as
\begin{equation}
\s_u \qp \s_v = \sum_{w\in S_n} 
                C_{u,v,w}\, \s_{\wnot w}\,,
\label{eq:su*sv}
\end{equation}
where $C_{u,v,w}\in\Z[q_1,\dots,q_{n-1}]$
and 
$\wnot=
\left(
\mbox{\scriptsize
$
\begin{array}{cccc}
1 & 2 & \cdots & n\\
n & n-1 & \cdots & 1
\end{array}
$
} \right)
$
is the longest permutation in~$S_n$.

The coefficient of $q_1^{d_1}\cdots q_{n-1}^{d_{n-1}}$
in the polynomial $C_{u,v,w}$ 
is the {\it Gromov-Witten invariant} 
$\left<\sigma_u,\sigma_v,\sigma_w\right>_{(d_1,\dots,d_{n-1})}$.
The Gromov-Witten invariants are defined geometrically as
numbers of certain rational curves in~$Fl_n$.  (See~\cite{FP}
or~\cite{FGP} for details.)
Let us summarize the main properties of these invariants.
It will be more convenient for us to work with the polynomials $C_{u,v,w}$.

{\it 
\begin{enumerate}
\item[{\bf 1.}] {\rm (Nonnegativity)} \ 
     All coefficients of the $C_{u,v,w}$ are nonnegative integers.
\item[{\bf 2.}] {\rm ($S_3$-symmetry)} \ 
     The polynomials $C_{u,v,w}$ are invariant with respect to permuting
     $u$, $v$, and~$w$.
\item[{\bf 3.}] {\rm (Degree condition)} \
     The polynomial $C_{u,v,w}$ is a homogeneous polynomial of degree 
     ${1\over 2} (\l(u)+\l(v)+\l(w)-{n(n-1)\over 2})$.
\item[{\bf 4.}] {\rm (Classical limit)} \
      The Schubert number $c_{u,v,w}$ is 
     the constant term of the polynomial $C_{u,v,w}$.
\item[{\bf 5.}] {\rm (Associativity)} \
      The operation~``\kern0.03cm$*$\kern-.05cm'' defined 
      by~{\rm(\ref{eq:su*sv})} via the polynomials $C_{u,v,w}$ is associative.
\end{enumerate}}

The first four properties are clear from geometric definitions.
It was conjectured in~\cite{FGP} that nonnegativity,
associativity, degree condition, and classical limit condition
uniquely determine the Gromov-Witten invariants.

The conditions~{\bf 3} and~{\bf 4} immediately imply the following 
statement.

\begin{proposition}
\label{prop:simple}
We have
$$
C_{u,v,w}=\left\{
\begin{array}{cl}

0        & \textrm{if }\l(u)+\l(v)+\l(w)<{n(n-1)\over 2},\\[.05in]
0        & \textrm{if }\l(u)+\l(v)+\l(w)-{n(n-1)\over 2} 
                                        \textrm{ is odd,}\\[.05in]
c_{u,v,w}& \textrm{if }\l(u)+\l(v)+\l(w)={n(n-1)\over 2},\\[.05in]
???& \textrm{overwise.}
\end{array}
\right.
$$
\end{proposition}

\medskip

In~\cite{FGP} we gave a method for calculation of the Gromov-Witten invariants.
Among several approaches presented in that paper, one is
based on the quantum analogue of Monk's formula.

For $1\leq i<j\leq n$, let $s_{ij}$ be the transposition in~$S_n$ that 
permutes $i$ and~$j$. 
Then $s_i = s_{i\,i+1}$ is an adjacent transposition.
Also, let $q_{ij}$ be a shorthand for the product $q_i q_{i+1} \cdots q_{j-1}$.

\begin{proposition}
\label{prop:qmonk} 
{\rm \cite[Theorem~1.3]{FGP} (quantum Monk's formula)} \ 
For $w\in S_n$ and $1\leq k<n$, 
the quantum product of the Schubert classes~$\s_{s_k}$ and~$\s_w$
is given by
\begin{equation}
\sigma_{s_k}\qp\sigma_w= \sum \s_{w s_{ab}} +  \sum q_{cd}\, \s_{w s_{cd}}\,,
\label{eq:qmonk}
\end{equation}
where the first sum is over all transpositions~$s_{ab}$ such that
$a\leq k<b$ and $\l(w s_{ab})=\l(w) +1$, and the second sum is over
all transpositions $s_{cd}$ such that $c\leq k<d$ and 
$\l(w s_{cd}) =\l(w)-\l(s_{cd}) = \l(w) - 2(d-c)+1$.
\end{proposition}

\begin{remark} 
\label{re:generate}
{\rm 
The two-dimensional Schubert classes~$\s_{s_k}$ 
generate the quantum cohomology ring.
Thus formula~(\ref{eq:qmonk}) uniquely determines the multiplicative structure
of~$\QH^*(Fl_n,\Z)$ and, therefore, the Gromov-Witten invariants.
}
\end{remark}

\section{Cyclic symmetry}
\label{se:main}

Let $o=(1,2,\dots,n)$ be the cyclic permutation in~$S_n$ given by 
$$
o(i)=i+1,\textrm{ for } i=1,\dots,n-1,\qquad o(n)=1.
$$

Recall that $q_{ij}=q_i q_{i+1} \cdots q_{j-1}$ for $i<j$.  We also 
define $q_{ij} = q_{ji}^{-1}$ for $i>j$ and $q_{ii}=1$.  

\begin{theorem}
\label{th:main} 
For any $u,v,w\in S_n$ we have
\begin{equation}
C_{u,v,w} = q_{ij}\, C_{u,o^{-1} v, o w}\,,
\label{eq:cuvw}
\end{equation}
where $i=v^{-1}(1)$ and $j=w^{-1}(n)$.
\end{theorem}

The $S_3$-invariance of the~$C_{u,v,w}$ under permuting $u$, $v$, and~$w$
implies a more general statement.

For $w\in S_n$ and $1\leq a\leq n$, define the following 
Laurent monomials in the $q_i$ 
$$
Q_{w,a} = \prod_{i\,:\,w(i)\geq n-a+1} q_{1i}\,,\qquad
Q_{w,-a} = \prod_{j\,:\,w(j)\leq a} (q_{1j})^{-1},
$$
and let $Q_{w,0}=1$.

\begin{corollary}
\label{cor:QQQ} 
 For any $u,v,w\in S_n$ and $-n\leq a,b,c\leq n$ such that 
$a+b+c=0$, we have
\begin{equation}
C_{u,v,w}= Q_{u,a} Q_{v,b} Q_{w,c}\, 
C_{o^a u,o^b v, o^c w}\,.
\label{eq:QQQ}
\end{equation}
\end{corollary}


In many cases Corollary~\ref{cor:QQQ} and Proposition~\ref{prop:simple} 
allow us to reduce the polynomials $C_{u,v,w}$ to the Schubert numbers 
$c_{u,v,w}$:

\begin{corollary}  For $u,v,w\in S_n$,
let us find a triple $-n\leq a,b,c\leq n$, $a+b+c=0$,
for which the expression
$$
\l_{a,b,c}=\l(o^a u)+\l(o^b v)+\l(o^c w)
$$
is as small as possible.  If $\l_{a,b,c}< {n(n-1)\over 2}$ then $C_{u,v,w}=0$.
If $\l_{a,b,c}={n(n-1)\over 2}$ then $C_{u,v,w} = Q_{u,a} Q_{v,b} Q_{w,c}\, 
c_{o^a u,o^b v, o^c w}$\,.
\end{corollary}

%

\begin{remark} 
{\rm 
(Reduction of Gromov-Witten invariants) \
The Gromov-Witten invariants have the following {\it stability property.}
If $u,v,w\in S_n$ are three permutations such that $u(n)=v(n)=n$ and
$w(n)=1$ then $C_{u,v,w}=C_{u',v',w'}$, where $u',v',w'\in S_{n-1}$
are permutations obtained from $u,v,w$ by removing the last entry 
(and subtracting 1 from all entries of $w$).

For a triple of permutation $u,v,w\in S_n$ such that $u(n)+v(n)+w(n)\equiv 1
\pmod n$, we can use the relation~(\ref{eq:QQQ}) to transform the triple
to the above case when we can use the stability property.
This shows that $1/n$ of all Gromov-Witten invariants for $Fl_n$ can
be reduced to the Gromov-Witten invariants of~$Fl_{n-1}$.
Analogously, we can reduce the problem to a lower level for a triple
of permutations $u,v,w\in S_n$ such that $u(1)+v(1)+w(1)\equiv 2\pmod n$.
}
\end{remark}

\begin{remark} 
{\rm 
(New rules for multiplication of Schubert classes) \
Suppose that a rule is know for the quantum multiplication 
of an arbitrary Schubert class by certain Schubert class $\s_u$.
Theorem~\ref{th:main} immediately produces a new rule for the 
quantum multiplication by $\s_{o^a u}$, where $a\in \Z$.
For example, we get for free a rule for $\s_{o^a}*\s_w$.
Quantum Monk's formula~(\ref{eq:qmonk}) can be
extended to a rule for $\s_{o^a s_k}*\s_w$.
More generally, 
quantum Pieri's formula~\cite[Corollary~4.3]{P} extends to
an explicit rule for $\s_{o^a u}*\s_w$, where $u$ is a permutation
of the form $u=s_k s_{k+1} \cdots s_{k+l}$ or 
$u=s_k s_{k-1} \cdots s_{k-l}$.
}
\end{remark}

\section{Twisted cyclic shift}
\label{se:proof}

Let $T_{ij}$, $1\leq i<j\leq n$, be the $\Z[q_1,\dots,q_{n-1}]$-linear 
operators that act
on the quantum cohomology ring~$\QH^*(Fl_n,\Z)$ by
\begin{equation}
T_{ij}\,:\, \s_w\longmapsto 
\left\{
\begin{array}{cc}
\s_{w s_{ij}}           & \textrm{if } \l(w s_{ij}) = \l(w) +1,\hfill    \\
q_{ij}\, \s_{w s_{ij}}  & \textrm{if } \l(w s_{ij}) = \l(w) -2(j-i)+1,\hfill \\
0                       & \textrm{otherwise.}\hfill
\end{array}
\right.
\label{eq:tij}
\end{equation}
Then quantum Monk's formula~(\ref{eq:qmonk}) can be written as:
\begin{equation}
\label{eq:qMonk-T}
\s_{s_k}\qp\s_w = \sum_{i\leq k<j} T_{ij}(\s_w).
\end{equation}


The operators~$T_{ij}$ satisfy certain simple quadratic relations.  
The formal algebra defined by these relations was studied in~\cite{FK} 
and~\cite{P}.

Let us also define the {\it twisted cyclic shift operator\/}~$O$
that acts on the quantum cohomology ring $\QH^*(Fl_n,\Z)$, 
linearly over $\Z[q_1,\dots,q_{n-1}]$, by
$$
O\,:\, \s_w \longmapsto  q^{(w)}\, \s_{o w}\,,
$$
where $q^{(w)}=q_{rn}$ with $r=w^{-1}(n)$.

\begin{proposition}
\label{prop:commute}
For any $1\leq i<j\leq n$, the operators $T_{ij}$ and~$O$ commute:
$$
      T_{ij}\,O=O\,T_{ij}\,.
$$
\end{proposition}

The following lemma clarifies the conditions in the right-hand side 
of~(\ref{eq:tij}).  Its proof is a straightforward observation.
\begin{lemma}
\label{le:ijk} 
Let $w\in S_n$ and $1\leq i<j\leq n$.  Then
\begin{enumerate}
\item
$\l(w\,s_{ij})=\l(w)+1$ if and only if for all $i\leq k\leq j$ we have 
$$w(k)\geq w(j)\geq w(i)\quad \textrm{or}\quad w(j)\geq w(i)\geq w(k)\,;$$
\item 
$\l(w\,s_{ij})=\l(w) - \l(s_{ij})= \l(w)-2(j-i)+1$ 
if and only if for all $i\leq k\leq j$ we have
$$w(i)\geq w(k)\geq w(j)\,.$$
\end{enumerate}
\end{lemma}

\proofof{Proposition~\ref{prop:commute}} 
The crucial observation is that, for fixed $i\leq k\leq j$, 
the set of permutations~$w$ such that
$$
w(k)\geq w(j)\geq w(i) \quad \textrm{or}\quad 
w(j)\geq w(i)\geq w(k) \quad \textrm{or}\quad
w(i)\geq w(k)\geq w(j)
$$
is invariant under the left multiplications of $w$ by the cycle~$o$.
This fact, together with Lemma~\ref{le:ijk}, implies that $(T_{ij}\,O)(\s_w)$
is nonzero if and only if $T_{ij}(\s_w)$ is nonzero.
Assume that $T_{ij}(\s_w)\ne 0$ and consider three cases:

I. Neither $w(i)$ nor $w(j)$ is equal to $n$.
Then either of the conditions in the right-hand side of~(\ref{eq:tij})
is satisfied for $w$ if and only if the same condition is satisfied for $ow$.
Also $q^{(w)}=q^{(w s_{ij})}$.  Thus $(T_{ij}\,O)(\s_w)=(O\,T_{ij})(\s_w)$.

II. We have $w(j)=n$.  Then $w(i)<w(j)$ and $ow(i)> ow(j)$.  
Thus $\l(w s_{ij}) = \l(w) + 1$ and $\l(o w s_{ij}) = \l(o w) - \l(s_{ij})$.
Thus  $T_{ij}(\s_w)= \s_{w s_{ij}}$ and 
$T_{ij}(\s_{ow})= q_{ij} \s_{o w s_{ij}}$.  Also we have $q^{(w)} = q_{jn}$
and $q^{(w s_{ij})} = q_{in}$.
Therefore, $(T_{ij}\,O)(\s_w)= q_{ij} q_{jn} \s_{o w s_{ij}}  =
q_{in} \s_{o w s_{ij}} = (O\, T_{ij})(\s_w)$.

III.  We have $w(i)=n$.  Then $w(i)>w(j)$ and $ow(i)< ow(j)$.  
Thus $\l(w s_{ij}) = \l(w) - \l(s_{ij})$ and $\l(o w s_{ij}) = \l(o w) +1$.
Thus  $T_{ij}(\s_w)= q_{ij} \s_{w s_{ij}}$ and 
$T_{ij}(\s_{ow})= \s_{o w s_{ij}}$.  Also we have $q^{(w)} = q_{in}$
and $q^{(w s_{ij})} = q_{jn}$.
Therefore, $(T_{ij}\,O)(\s_w)= q_{in} \s_{o w s_{ij}}  =
q_{ij} q_{jn} \s_{o w s_{ij}} = (O\, T_{ij})(\s_w)$.
\endproof

\begin{corollary} 
\label{cor:commute}
For any~$w\in S_n$, the operator of quantum multiplication by the 
Schubert class~$\s_{w}$ commutes with the operator~$O$.  
\end{corollary}
\proof
Proposition~\ref{prop:commute} and quantum Monk's formula~(\ref{eq:qMonk-T})
imply that the operator of quantum multiplication
by a two-dimensional Schubert class~$\s_{s_k}$ commutes with 
the twisted cyclic shift operator~$O$.
By Remark~\ref{re:generate}, for any~$w\in S_n$, 
the operator of quantum multiplication by~$\s_{w}$ commutes with~$O$.  
\endproof

This also proves Theorem~\ref{th:main}, because
it is equivalent to Corollary~\ref{cor:commute}.

\section{Transition graph}
\label{se:graph} 

The {\it Bruhat order\/} $Br_n$ is the partial order on the set of 
all permutations in $S_n$ given by the following covering relation:
$u\to w$ if $w=u\,s_{ab}$ and $\l(w)=\l(u)+1$.
In other words, $u\to w$ if $\s_w$ appear in the expansion
of $\s_{s_k}\cdot \s_u$ for some $1\leq k< n$ (the product in the usual 
cohomology ring).

The analogue of the Bruhat order for the quantum cohomology ring is the
following transition graph.  The {\it transition graph\/} $Tr_n$ is the
directed graph on the set of permutations in $S_n$.  Two permutations 
are connected by an edge $u\to w$ in $Tr_n$ if $w=u\,s_{ab}$ and 
either $\l(w)=\l(u)+1$ or $\l(w)=\l(u)-\l(s_{ab})$.  We will
label the edge $u\to u\, s_{ab}$ by the pair $(a,b)$.  Equivalently,
two permutations are connected by the edge $u\to w$ in $Tr_n$ 
whenever $\s_w$ appear in the expansion of the quantum product 
$\s_{s_k} * \s_u$ for some $1\leq k< n$.

Proposition~\ref{prop:commute} implies the cyclic symmetry
of the transition graph:

\begin{corollary}
The transition graph $Tr_n$ is invariant under the cyclic shift:
$w\mapsto o\,w$, for $w\in S_n$.
\end{corollary}

\begin{figure}[ht]
$$
\xymatrix{
&& *+[F-:<3pt>]{321} &&\\
*+[F-:<3pt>]{231}  \ar[urr]|{12} &&&& 
*+[F-:<3pt>]{312}  \ar[ull]|{23}  \\
\\
*+[F-:<3pt>]{213}  \ar[uu]|{23}  \ar[uurrrr]|(.3){13}  &&&& 
*+[F-:<3pt>]{132}  \ar[uu]|{12}  \ar|(.3){13}[uullll]   \\
&& *+[F-:<3pt>]{123} \ar[ull]|{12}  \ar[urr]|{23}    &&
}
$$
\caption{Bruhat order $Br_3$.}
\end{figure}
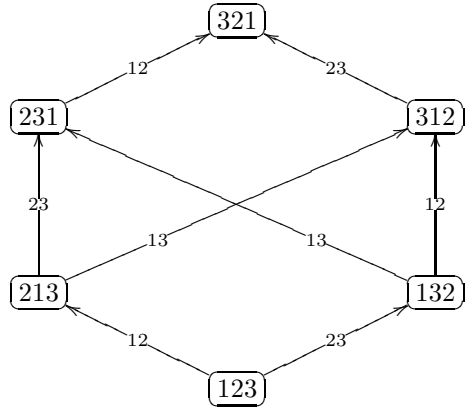
\bigskip
\bigskip

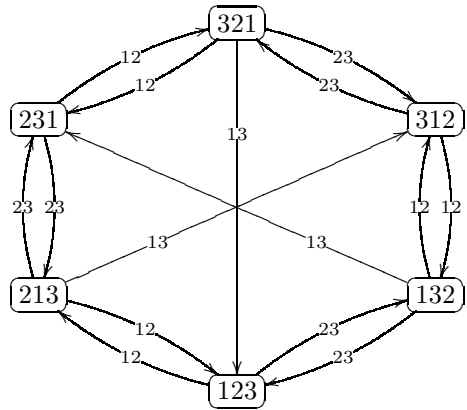
\begin{figure}[ht]
$$
\xymatrix{
&& *+[F-:<3pt>]{321}  \ar[dddd]|(.3){13}  \ar@/^/[dll]|{12}  \ar@/^/[drr]|{23}  &&\\
*+[F-:<3pt>]{231}  \ar@/^/[urr]|{12}  \ar@/^/[dd]|{23} &&&& 
*+[F-:<3pt>]{312}  \ar@/^/[ull]|{23}  \ar@/^/[dd]|{12}   \\
\\
*+[F-:<3pt>]{213}  \ar@/^/[uu]|{23}  \ar[uurrrr]|(.3){13}  \ar@/^/[drr]|{12} &&&& 
*+[F-:<3pt>]{132}  \ar@/^/[uu]|{12}  \ar|(.3){13}[uullll]  \ar@/^/[dll]|{23}  \\
&& *+[F-:<3pt>]{123} \ar@/^/[ull]|{12}  \ar@/^/[urr]|{23}    &&}
$$
\caption{Transition graph $Tr_3$.}
\end{figure}

\bigskip
\bigskip

Figures~1 and~2 
show the Bruhat order $Br_3$ and the transition graph $Tr_3$. 
The transition graph $Tr_3$ is obtained by adding several new edges
to $Br_3$, which makes the picture symmetric with respect
to the cyclic group $\Z/3\Z$.  The generator $o$ of the cyclic group
rotates the graph $Tr_3$ by $180^{\mathrm{o}}$ clockwise.

\end{document}